# ON ONE INVERSE PROBLEM

## Aliyev N.A., Mustafayeva Y.Y., Efendiyev R.F.

Baku State University,

23, Z.Khalilov str., AZ 1148, Baku, Azerbaijan

aliyev.nihan@mail.ru, helenmust@rambler.ru, rakibefendiev@bsu.edu.az

**Abstract**. The stated paper is dedicated to one of the inverse problems of spectral theory. It is necessary to define matrix (constant) coefficients of some quadratic pencil, if the eigenvalues of this pencil are known.

Furthermore, it is known that these eigenvalues are complex valued and besides that they contain their conjugate ones.

As is known, it is a problem of one company. We shall approach this problem on some restrictions.

**Introduction**.

We are starting from spectral theory of finite-dimensional space, i.e. from spectral theory of quadratic symmetric matrices [1]. So, let

$$A = \begin{pmatrix} a_{11} & a_{12} & ... & a_{1n} \\ a_{21} & a_{22} & ... & a_{2n} \\ ... & ... & ... & ... \\ a_{n1} & a_{n2} & ... & a_{nn} \end{pmatrix}$$

be a quadratic matrix of n-th order, $a_{ij} \in R$ and $a_{ij} = a_{ij}$, $i \neq j$, i.e. $A' = A$ where $A'$ is a transposed matrix. Let us define eigenvalues and eigenvectors of the matrix $A$.

Let us consider a system of linear algebraic equations

$$Ax = \lambda x$$

or

$$(A - \lambda I)x = 0 \qquad (1)$$

where $I$ is a unit matrix of n-th order. As is known, for the existence of nontrivial solution of system (1) the parameter $\lambda$ must satisfy the following equation (2):

$$\det(A - \lambda I) = |A - \lambda I| = 0 \qquad (2)$$

as $|A - \lambda I|$ is a polynomial of power $n$ with respect to $\lambda$. Let us assume that the roots of equation (2) are real-valued and different, i.e.

$$\lambda_k \in R,\ k = \overline{1,n}\ ;\ \lambda_k \neq \lambda_m,\ k \neq m,\ k.m = \overline{1,n}. \qquad (3)$$

Now we shall define the eigenvectors of the matrix A corresponding the eigenvalues $\lambda_k$. For this we replace $\lambda$ with $\lambda_k$ in the system of linear algebraic equations (1) and solving we get

$$X^{(k)} = (A_{g_k 1}, A_{g_k 2}, ..., A_{g_k n})',\ k = \overline{1,n}, \qquad (4)$$

where an element $A_{g_k j}$, $j = \overline{1,n}$, of the column $X^{(k)}$ is a cofactor of the element at the intersection of the $g_k$-th row and $j$-th column of the determinant $|A - \lambda I|$ and $g_k$ is the number of such a row for which the cofactors of all its elements don't turn into a zero.

Thus, the eigenvector of system (1), or the matrix A, corresponding to the eigenvalue $\lambda_k$ is built in the form of (4). Indeed, if in system (1) to replace the parameter $\lambda$ with $\lambda_k$ and the vector-column $X$ with the eigenvector $X^{(k)}$ in the form (4) then we have

$$\begin{pmatrix} a_{11} - \lambda_k & a_{12} & \cdots & a_{1n} \\ a_{21} & a_{22} - \lambda_k & \cdots & a_{2n} \\ \cdots & \cdots & \cdots & \cdots \\ a_{n1} & a_{n2} & \cdots & a_{nn} - \lambda_k \end{pmatrix} \begin{pmatrix} A_{g_k 1} \\ A_{g_k 2} \\ \cdots \\ A_{g_k n} \end{pmatrix} = 0. \qquad (5)$$

The validity of (5) follows from that if to select $g_k$-th row of the matrix then the product of this row with the column gives us the determinant $|A - \lambda I|$ which turns into zero as $\lambda_k$ is a zero of the determinant (2). And what relates to the product of each left row by the vector-column then their turning into zero follows from the properties of determinants. Now let us reduce the vector-column (4) to the normal form, i.e.

$$x^{(k)} = \frac{X^{(k)}}{\|X^{(k)}\|} = \left( \frac{A_{g_k 1}}{\|X^{(k)}\|}, \frac{A_{g_k 2}}{\|X^{(k)}\|}, ..., \frac{A_{g_k n}}{\|X^{(k)}\|} \right)',\ k = \overline{1,n}, \qquad (6)$$

As is known, if to multiply a row by a column of the same size then we'll get a scalar (it is a scalar product of vectors).

If to multiply a column by a row then we'll get a matrix, i.e.

$$(a_1, a_2, ..., a_n)'(b_1, b_2, ..., b_m) =$$

$$= \begin{pmatrix} a_1 b_1 & a_1 b_2 & ... & a_1 b_m \\ a_2 b_1 & a_2 b_2 & ... & a_2 b_m \\ ... & ... & ... & ... \\ a_n b_1 & a_n b_2 & ... & a_n b_m \end{pmatrix}$$

having size $n \times m$.

Then on the basis of (6) we build the following matrix-projectors (operators):

$$P_i = x^{(i)} \cdot x^{(i)'}, \quad i = \overline{1,n}, \qquad (7)$$

having the form of a square matrix of order $n$. Taking into account that the eigen vectors (6) are normal we have:

$$x^{(i)'} \cdot x^{(j)} = \delta_{ij} = \begin{cases} 0, & i \neq j, \\ 1, & i = j. \end{cases} \qquad (8)$$

So,

$$P_i \cdot P_i = x^{(i)} \cdot x^{(i)'} \cdot x^{(i)} \cdot x^{(i)'} = x^{(i)} \cdot \left[ x^{(i)'} \cdot x^{(i)} \right] \cdot x^{(i)'} = x^{(i)} \cdot x^{(i)'} = P_i,$$

i.e.

$$P_i^2 = P_i$$

and

$$P_i \cdot P_j = x^{(i)} \cdot x^{(i)'} \cdot x^{(j)} \cdot x^{(j)'} = x^{(i)} \cdot \left[ x^{(i)'} \cdot x^{(j)} \right] \cdot x^{(j)'} = 0, \quad i \neq j. \qquad (9)$$

Let $y$ is an arbitrary vector from the domain of definition of problem (1). Then expanding it by basis (6) we have:

$$y = \sum_{i=1}^{n} y_i \, x^{(i)}, \tag{10}$$

where $y_i$ are constant coefficients.

Taking into account that

$$\sum_{i=1}^{n} y_i \, x^{(i)} = \sum_{i=1}^{n} x^{(i)} y_i \, ,$$

the constants $y_i$ can be represented in the form:

$$y_i = x^{(i)'} y, \ i = \overline{1, n}. \tag{11}$$

Then on the basis of the equality

$$P = \sum_{i=1}^{n} P_i \, , \tag{12}$$

we have

$$P y = \sum_{i=1}^{n} P_i \, y = \sum_{i=1}^{n} x^{(i)} \cdot x^{(i)'} y = \sum_{i=1}^{n} x^{(i)} y_i = \sum_{i=1}^{n} y_i x^{(i)} = y,$$

i.e.

$$P = \sum_{i=1}^{n} P_i = I \, , \tag{13}$$

where $I$ is a unit matrix of the n-th order.

Now it is easy to see that

$$\sum_{i=1}^{n} \lambda_i \, P_i = \sum_{i=1}^{n} \lambda_i \, x^{(i)} \cdot x^{(i)'} = \sum_{i=1}^{n} \left( \lambda_i \, x^{(i)} \right) \cdot x^{(i)'} . \tag{14}$$

As $\lambda_i$ is an eigenvalue of the matrix $A$ and $x^{(i)}$ is the corresponding eigenvector then from (1) we get

$$\lambda_i \, x^{(i)} = A x^{(i)}, \ i = \overline{1, n}.$$

Then from (14) we have

$$\sum_{i=1}^{n}\lambda_i P_i = \sum_{i=1}^{n} A x^{(i)} \cdot x^{(i)'} = A\sum_{i=1}^{n} x^{(i)} \cdot x^{(i)'} = A\sum_{i=1}^{n} P_i = A \cdot I = A.$$

Thus, we have obtained a spectral expansion for a real-valued symmetric matrix, i.e. spectral expansion of a self-adjoint operator:

$$A = \sum_{i=1}^{n} \lambda_i P_i . \tag{15}$$

Consequently, we can easily see that for any natural number $k$ we have

$$A^k = \sum_{i=1}^{n} \lambda_i^k P_i . \tag{16}$$

Therefore, we have proved the following statement:

**Theorem**. *For any real-valued finite-dimensional symmetric matrix A with different real-valued eigenvalues there hold true relationships* (16).

Now we start a new stage of this work. Let a quadratic pencil have the form:

$$P(\lambda) = \lambda^2 M + \lambda C + K , \tag{17}$$

where

$$C = C(\alpha) = C_0 + \sum_{i=1}^{n} \alpha_i C_i , \tag{18}$$

$$K = K(\beta) = K_0 + \sum_{i=1}^{n} \beta_i K_i , \tag{19}$$

$$\alpha = (\alpha_1, \alpha_2, ..., \alpha_n) \in R^n , \quad \beta = (\beta_1, \beta_2, ..., \beta_n) \in R^n ,$$

$M, C_i, K_i$ are real-valued symmetric matrices of $n$-th order. $\{\lambda_1, \lambda_2, ..., \lambda_{2n}\}$ are complex-valued eigenvalues of a quadratic pencil (17) where each eigenvalue is included with its conjugate, i.e.

$$\lambda_1(\alpha,\beta), \overline{\lambda_1}(\alpha,\beta), \lambda_2(\alpha,\beta), \overline{\lambda_2}(\alpha,\beta), ,..., \lambda_n(\alpha,\beta), \overline{\lambda_n}(\alpha,\beta) \tag{20}$$

are eigenvalues of the quadratic pencil (17). Thus, the set (20) is a set of the roots of the equation

$$\det P(\lambda) = \det((\lambda^2 M + \lambda C + K) = 0 . \tag{21}$$

we suggest that $\lambda_k(\alpha, \beta), k = \overline{1,n}$, are distinct, i.e.

$$\lambda_k(\alpha, \beta) \neq \lambda_m(\alpha, \beta), \text{ if } m \neq k, m, k = \overline{1,n}. \tag{22}$$

Then suggesting that

$$C_0 = 0, K_0 = 0,$$

$\alpha_i, i = \overline{1,n}, (\alpha_i \in R)$, can be taken as eigenvalues of the matrix $C(\alpha)$, and matrices $C_i$ as the projectors of the matrix $C(\alpha)$ correspond to the eigenvalue $\alpha_i$. Similarly, $\beta_i, i = \overline{1,n}, (\beta_i \in R)$, can be taken as eigenvalues of the matrix $K(\beta)$, and matrices $K_i$ as the projectors of the matrix $K(\beta)$ correspond to the eigenvalue $\beta_i$.

If $\alpha_i$ and $\beta_i, i = \overline{1,n}$, are known then $\lambda_k, \overline{\lambda_k}, k = \overline{1,n}$, are defined from the equation (21).

Now we suggest that set (20) is known and condition (22) holds true. The elements of set (20) are eigenvalues of a quadratic pencil (17).

Let matrix M be non-degenerated, i.e.

$$\det M \neq 0. \tag{23}$$

Then from (17) we obtain that

$$\tilde{P}(\lambda) = M^{-1}P(\lambda) = \lambda^2 I + \lambda \tilde{C} + \tilde{K} \tag{24}$$

where

$$\tilde{C} = M^{-1}C, \tilde{K} = M^{-1}K. \tag{25}$$

Suggesting that $\tilde{C}$ and $\tilde{K}$ are commutative matrices, i.e.

$$\tilde{C}\tilde{K} = \tilde{K}\tilde{C}. \tag{26}$$

Then their eigenvectors coincide [2]-[4], so coincide their projectors and the orthogonal matrix which reduces them to diagonal form.

Let $\tilde{\alpha}_i \in R$ $i = \overline{1,n}$, be eigenvalues of the matrix $\tilde{C}$ and $\tilde{\beta}_i \in R$, $i = \overline{1,n}$, be eigenvalues of the matrix $\tilde{K}$. Let us designate the normed eigenvectors corresponding to these eigenvalues as $\tilde{z}_i$. Then from (26) we have

$$\tilde{C}\tilde{K}\tilde{z}_i = \tilde{K}\tilde{C}\tilde{z}_i. \tag{27}$$

Taking into account that

$$\tilde{K}\tilde{z}_i = \tilde{\beta}_i \tilde{z}_i \text{ and } \tilde{C}\tilde{z}_i = \tilde{\alpha}_i \tilde{z}_i,$$

from (27) we obtain:

$$\tilde{C}\tilde{\beta}_i \tilde{z}_i = \tilde{K}\tilde{\alpha}_i \tilde{z}_i$$

or

$$\tilde{\beta}_i \tilde{C}\tilde{z}_i = \tilde{\alpha}_i \tilde{K}\tilde{z}_i$$

or we come to the identity

$$\tilde{\beta}_i \tilde{\alpha}_i \tilde{z}_i = \tilde{\alpha}_i \tilde{\beta}_i \tilde{z}_i.$$

Then designating the projectors of the matrices $\tilde{C}$ and $\tilde{K}$ as

$$\tilde{P}_1, \tilde{P}_2, ..., \tilde{P}_n \tag{28}$$

we get

$$\tilde{C} = \sum_{i=1}^{n} \tilde{\alpha}_i \tilde{P}_i, \tag{29}$$

$$\tilde{K} = \sum_{i=1}^{n} \tilde{\beta}_i \tilde{P}_i, \tag{30}$$

and

$$I = \sum_{i=1}^{n} \tilde{P}_i. \tag{31}$$

Taking into account (29)-(31), we find from (24) the following:

$$\tilde{P}(\lambda) = \sum_{i=1}^{n} (\lambda^2 + \tilde{\alpha}_i \lambda + \tilde{\beta}_i)\tilde{P}_i. \tag{32}$$

From the obtained expression (32) it follows that the eigenvalues of the quadratic pencil are obtained from the equation:

$$\lambda^2 + \tilde{\alpha}_i \lambda + \tilde{\beta}_i = 0, \ i = \overline{1,n}, \qquad (33)$$

for which $\tilde{\alpha}_i \in R, \ \tilde{\beta}_i \in R, \ i = \overline{1,n}$, and $\lambda_k \in C, \ k = \overline{1,n}, \ -\tilde{\alpha}_i = \lambda_i + \overline{\lambda}_i, \ \tilde{\beta}_i = \lambda_i \cdot \overline{\lambda}_i, \ i = \overline{1,n}$.

By these relationships there is established the connection between the values of $\tilde{\alpha}_i$, $\tilde{\beta}_i$ and $\lambda_i$.

Now we'll consider the inverse problem.

Let us be given a set (20) satisfying conditions (22).

Then we build real-valued numbers $\tilde{\alpha}_i$ and $\tilde{\beta}_i$, $i = \overline{1,n}$, as follows:

$$\tilde{\alpha}_i = -(\lambda_i + \overline{\lambda}_i), \ \tilde{\beta}_i = \lambda_i \cdot \overline{\lambda}_i, \ i = \overline{1,n}.$$

After that the numbers $\tilde{\alpha}_i$ and $\tilde{\beta}_i$, $i = \overline{1,n}$, are defined, the matrices $\tilde{C}$ and $\tilde{K}$ will be defined by means of (29) and (30).

For this it is necessary to find projectors. Taking into account that each projector is one-dimensional (i.e. the rank of the projector is unit), they can be taken in the form

$$\tilde{P}_k = \begin{pmatrix} 0 & . & . & . & . & . & 0 \\ 0 & . & . & . & . & . & 0 \\ . & . & . & . & . & . & . \\ 0 & . & 0 & 1 & 0 & . & 0 \\ . & . & . & . & . & . & . \\ 0 & 0 & . & . & . & . & 0 \end{pmatrix} \leftarrow k-th \ row$$

$$\uparrow$$
$$k - th \ column$$

**Remark**. Matrices $\tilde{C}$ and $\tilde{K}$ can be obtained from the diagonal matrices

$$\begin{pmatrix} \tilde{\alpha}_1 & .. & 0 \\ .. & . & .. \\ 0 & ... & \tilde{\alpha}_n \end{pmatrix} \text{ and } \begin{pmatrix} \tilde{\beta}_1 & .. & 0 \\ .. & . & .. \\ 0 & ... & \tilde{\beta}_n \end{pmatrix}$$

by means of any orthogonal matrix.